\documentclass[draft]{publmathd}
\usepackage{amsmath,amsfonts,amssymb}

\theoremstyle{plain}
\newtheorem{lemma}{Lemma}
\newtheorem{Theorem}{Theorem}
\newtheorem{Corollary}{Corollary}

\DeclareMathOperator{\charac}{char}
\DeclareMathOperator{\expon}{exp}
\def\gp#1{\langle #1 \rangle}
\def\m1{^{-1}}

\author{VICTOR BOVDI, TIBOR JUH\'ASZ  and  ERNESTO SPINELLI}

\address{
Institute of Mathematics, University of Debrecen, \\
H-4010 Debrecen, P.O.B. 12, Hungary\\
Institute of Mathematics and Informatics,\quad College of Ny\'\i
regyh\'aza\\ S\'ost\'oi \'ut 31/b, H-4410 Ny\'\i regyh\'aza,
Hungary} \email{vbovdi@math.klte.hu }

\address{
Institute of Mathematics, University of Debrecen, \\
H-4010 Debrecen, P.O.B. 12, Hungary} \email{juhaszti@math.klte.hu
}

\address{
Dipartimento di Matematica "E. De Giorgi",\\
Universit\`{a} degli Studi di Lecce\\
Via Provinciale Lecce-Arnesano, \quad 73100-LECCE, Italy}
\email{spinelli@ilenic.unile.it\\}

\title[Modular  group algebras with  almost maximal Lie nilpotency indices. I]
      {Modular  group algebras with \\ almost maximal Lie nilpotency indices. I}
\thanks{The research was supported by OTKA  No.T 037202, No.T 038059 and Italian
National Research Project "Group Theory and Application"}



\submitted{March 8, 2005}

\begin{document}

\begin{abstract}
Let $K$ be a field of positive characteristic $p$ and $KG$ the
group algebra of a group $G$. It is known that, if $KG$ is Lie
nilpotent, then its upper (or lower) Lie nilpotency index is at
most $|G'|+1$, where $|G'|$ is the order of the commutator
subgroup.  The authors have previously determined the groups $G$
for which this index is maximal and here they determine the $G$
for which it is `almost maximal', that is the next highest
possible value, namely  $|G'|-p+2$.
\end{abstract}
\maketitle
%\section{Introduction and results}
Let $R$ be an associative algebra with identity.  The algebra $R$
can be regarded as a Lie  algebra, called the associated Lie
algebra of $R$, via the  Lie  commutator $[x,y]=xy-yx$, for every
$x,y\in R$. Set $[x_1,\ldots,x_{n}]=[[x_1,\ldots,x_{n-1}],x_{n}],$
where $x_1,\ldots,x_{n}\in R$. The \emph{$n$-th lower Lie power}
$R^{[n]}$ of $R$ is the associative ideal generated by all the Lie
commutators $[x_1,\ldots,x_n]$, where  $R^{[1]}=R$ and
$x_1,\ldots,x_n \in R$. By induction, we define the \emph{$n$-th
upper  Lie power}  $R^{(n)}$ of $R$ as the associative ideal
generated by all the Lie commutators $[x,y]$, where $R^{(1)}=R$
and $x\in R^{(n-1)}$, $y\in R$.

The algebra  $R$ is called  \emph{Lie nilpotent} (respectively
\emph{upper Lie nilpotent}) if there exists $m$ such that
$R^{[m]}=0$\quad  ($R^{(m)}=0$). The algebra  $R$ is called
\emph{Lie hypercentral} if for each sequence $\{a_i\}$ of elements
of $R$ there exists some $n$ such that $[a_1,\ldots,a_n]=0$. The
minimal integers $m,n$ such that $R^{[m]}=0$ and $R^{(n)}=0$  are
called \emph{ the Lie nilpotency index} and \emph{ the upper Lie
nilpotency  index} of $R$ and they are denoted by $t_{L}(R)$ and
$t^{L}(R)$, respectively.

Let $U(KG)$ be the  group of units  of a group algebra $KG$. For
the noncommutative modular group algebra $KG$  the following
Theorem due to A.A.~Bovdi, I.I.~Khripta, I.B.S.~Passi,
D.S.~Passman and etc. (see \cite{BKh,K}) is well known: The
following statements are equivalent:
(a)\; $KG$ is Lie nilpotent;
(b)\; $KG$ is Lie hypercentral;
(c)\; $KG$ is strongly Lie nilpotent;
(d)\; $U(KG)$ is nilpotent;
(e)\; $\charac(K)=p>0$, $G$ is nilpotent
and its commutator subgroup $G^{\prime}$ is a finite $p$-group.

It is well known (see \cite{P}) that,  if $KG$ is Lie nilpotent,
then
$$
t_{L}(KG)\leq t^{L}(KG)\leq \vert G^{\prime}\vert +1.
$$
Moreover, according to \cite{BP}, if $\charac(K)>3$, then
$t_{L}(KG)= t^{L}(KG)$. But the question of when $t_{L}(KG)=
t^{L}(KG)$ for $\charac(K)=2,3$ is in general still open.
Several  important results on this topic were obtained in
\cite{BK}.

We say that  a Lie nilpotent group algebra $KG$  has \emph{upper
almost maximal}   Lie nilpotency index,\quad   if $t^{L}(KG)=\vert
G^{\prime}\vert-p+2$.

A.Shalev in \cite{S2} began study the question when the Lie
nilpotent group algebras $KG$ has  the maximal upper Lie
nilpotency index. In \cite{BS} was given the complete description
of such  Lie nilpotent group algebras. Using  these results  we
proof the following
\begin{Theorem}\label{Th1}
Let $KG$ be a Lie nilpotent group algebra  over a field $K$ of
positive characteristic  $p$. Then $KG$ has  upper almost maximal
Lie nilpotency index  if and only if one of the following
conditions holds:
\begin{enumerate}
\item[(i)] $p=2$,\quad  $cl(G)=2$ and $\gamma_2(G)$  is noncyclic of order $4$;
\item[(ii)] $p=2$,\quad  $cl(G)=4$,  $\gamma_2(G)=C_4\times C_2$ and $\gamma_3(G)=C_2\times C_2$;
\item[(iii)] $p=2$,\quad  $cl(G)=4$,  $\gamma_2(G)$ is elementary abelian of order $8$;
\item[(iv)] $p=3$,\quad  $cl(G)=3$ and $\gamma_2(G)$  is  elementary abelian of order $9$.
\end{enumerate}
\end{Theorem}

As a consequence, we obtain that the Theorem 3.9 of \cite{S2} can
not be extent  for $p=2$ and $p=3$:
\begin{Corollary} Let $K$ be a field with
$\charac(K)=p>0$ and $G$  a nilpotent group such that
$|G^{\prime}|=p^n$.
\begin{enumerate}
\item[(i)] If $p=2$ and $t_L(KG)<2^n+1$,  then $t_L(KG)\leq 2^{n}$.
\item[(i)] If $p=3$ and $t_L(KG)<3^n+1$,  then $t_L(KG)\leq 3^{n}-1$.
\end{enumerate}
\end{Corollary}

The authors would like to thank for Prof. L.G.~Kov\'acs for his
valuable comment and suggestions for clarifying the exposition.

%\section{Preliminaries}
We use the  standard   notation:\quad $C_n$ is the cyclic group of
order $n$;\; $\zeta(G)$ is the center of the group $G$,\;
$(g,h)=g^{-1}h^{-1}gh=g^{-1}g^h$\; ($g,h\in G$); \newline
$\mathrm{Q}_{8}=\gp{\; a,b\; \mid\;  a^{4}=1, \; b^2=a^{2}, \;
a^b=a^{-1}\;}$ the quaternion group of order $8$;
\newline
$\mathrm{D}_{8}=\gp{\; a,b\; \mid \; a^{4}=b^2=1, \; a^b=a^{-1}\;
}$ the dihedral group of order $8$;
\newline
$\gamma_{i}(G)$\quad is the $i$-{th} term of the lower central
series of $G$, i.e.
$$
\gamma_{1}(G)=G,\quad \quad
\gamma_{i+1}(G)=\big(\gamma_{i}(G),G\big)\quad \quad (i\geq 1).
$$
Let $K$ be a field of positive  characteristic $char(K)=p$ and $G$
a group. We consider  a sequence of subgroups of $G$, setting
$$
\mathfrak D_{(m)}(G)=G\cap ( 1+KG^{(m)}),\qquad\quad (m\geq 1).
$$
The subgroup $\mathfrak D_{(m)}(G)$ is  called the $m$-th   \emph{
Lie dimension subgroup} of $KG$. It is possible to describe  the
$\mathfrak D_{(m)}(G)$'s in terms of the lower central series of
$G$ in the following manner (see Theorem 2.8 of \cite{P}, p.48):
\begin{equation}\label{e:1}
\frak D_{(m+1)}(G)=\left \{
\begin{array}{ll}
G & \quad \text{if} \quad  m=0;\\
G^{\prime}& \quad \text{if}  \quad     m=1;\\
{\big(\frak D_{(m)}(G),G\big)(\frak D_{(\lceil {\frac{m}{p}}\rceil
+1)}(G))^p}& \quad \text{if}  \quad  m\geq 2,
\end{array} \right.
\end{equation}
where $\lceil {\frac{m}{p}}\rceil $ is the upper integer part of
${\frac{m}{p}}$.

By \cite{P} (see p.46) there exists an explicit expression for
$\frak D_{(m+1)}(G)$:
\begin{equation}\label{e:2}
\frak D_{(m+1)}(G)=\prod_{(j-1)p^i\geq m}\gamma_j(G)^{p^i}.
\end{equation}

Put $p^{d_{(k)}}=[\frak D_{(k)}(G):\frak D_{(k+1)}(G)]$, where
$k\geq 1$. If $KG$ is Lie nilpotent and $|G^{\prime}|=p^n$, then
according to Jennings' theory \cite{S3} for the Lie dimension
subgroups, we get
\begin{equation}\label{e:3}
\displaystyle t^{L}(KG)=2+(p-1) \sum_{m\geq 1} md_{(m+1)},
\end{equation}
and it is easy to check that
\begin{equation}\label{e:4}
\sum_{m\geq 2} d_{(m)}=n.
\end{equation}
We use the following results by Shalev (see Corollary 4.5 and
Corollary 4.6 of \cite{S1} and Theorem 3.9 of \cite{S2}):
\begin{lemma}\label{sha} Let $K$ be a field with
$\charac(K)=p>0$ and $G$  a nilpotent group such that
$|G^{\prime}|=p^n$ and  $\expon(G^{\prime})=p^l$.
\begin{enumerate}
\item[(i)]  If $d_{(m+1)}=0$ and $m$ is a power of $p$, then
$\frak D_{(m+1)}(G)=\gp{1}$.
\item[(ii)] If $d_{(m+1)}=0$ and $p^{l-1}$ divides
$m$, then $\frak D_{(m+1)}(G)=\gp{1}$.
\item[(iii)] If $p\geq 5$ and $t_L(KG)<p^n+1$,  then $t_L(KG)\leq
p^{n-1}+2p-1$.
\end{enumerate}
\end{lemma}
Fist of all we begin by proving the following:

\begin{lemma}\label{L:2}
Let $K$ be a field with $\charac(K)=p>0$ and $G$ a nilpotent group
such that  $\vert G^{\prime} \vert=p^{n}$. Then \quad
$t^{L}(KG)=\vert G^{\prime} \vert -p+2$\quad if and only if one of
the following conditions holds:
\begin{itemize}
\item [(i)]  $p=2$, \quad $n=2$\quad  and \quad $d_{(2)}=2$;
\item [(ii)]  $p=2$,\quad  $n>2$, \quad  $d_{(2^{i}+1)}=d_{(2^{n-1})}=1$ \quad and \quad $d_{(j)}=0$,
\newline
where \quad $0\leq i\leq n-2$, \quad $j\neq 2^i +1$, \quad $j\neq
2^{n-1}$\quad and $j>1$;
\item [(iii)] $p=3$, \quad $n=2$\quad  and \quad $d_{(2)}=1$\; $d_{(3)}=1$;
\item [(iv)] $p=3$, \quad  $n\geq 2$, \quad  $d_{(3^{i}+1)}=d_{(3^{n-1})}=1$ \quad and \quad $d_{(j)}=0$,
\newline
where \quad $0\leq i\leq n-2$, \quad $j\neq 3^i +1$, \quad $j\neq
3^{n-1}$\quad and $j>1$.
\end{itemize}
\end{lemma}
\begin{proof} If $n=1$ then according to Theorem 1 and Corollary 1
of \cite{BS}, $t_L(KG)=t^L(KG)=\vert G^{\prime} \vert+1$. So let
$p=n=2$  and $t^{L}(KG)=4$. By  (\ref{e:3}) it follows at once
that $d_{(2)}=2$. The converse is trivial.

Analogue, if $p=3$ and $n=2$ we get the statement (iv) of our
lemma.

Now suppose that either (ii) or (iv) of lemma  holds. Then, by
(\ref{e:3}), we get $t^{L}(KG)=\vert G^{\prime}\vert-p+2$, where
$p=2,3$. In order to prove the other implication, we preliminarily
state that $p\leq 3$.

Indeed, assume that $p\geq 5$ and $G^{\prime}$ is not cyclic. By
\cite{BP} and by (iii) of Lemma \ref{sha} we have that \;
$t^{L}(KG)=t_{L}(KG)\leq p^{n-1}+2p-1$. But $
p^{n}-p+2>p^{n-1}+2p-1$, because \quad $(p^{n-1}-3)(p-1)>0$.

Therefore we can assume that $n>2$ and either  $p=3$ or $p=2$.
First, we shall show that \quad $d_{(p^{i}+1)}>0$\quad for\quad
$0\leq i\leq n-2$.

Suppose there exists \quad  $0\leq s\leq n-2$ such that
$d_{(p^{s}+1)}=0$.  From (\ref{e:1}) it follows at once that
$s\neq 0$ and by (i) of Lemma \ref{sha} we have that $\frak
D_{(p^{s}+1)}(G)=\gp{1}$ and so $d_{(r)}=0$ for every $r\geq
p^{s}+1$. Moreover, if  $d_{({q}+1)}\not=0$, then $q\leq p^s$.
According to (\ref{e:3}) it follows that
\begin{align*}
t^{L}(KG)&=2+(p-1)(\sum_{i=0}^{s-1}{p^{i}}+\sum_{i=0}^{s-1}{p^{i}}(d_{({p^{i}}+1)}-1)+ \sum_{q\not=p^i}qd_{({q}+1)})\\
&<2+(p-1)(\sum_{i=0}^{s-1}{p^{i}}+(\sum_{i=0}^{s-1}(d_{({p^{i}}+1)}-1)+ \sum_{q\not=p^i}d_{({q}+1)})p^{s})\\
&=2+(p-1)(\sum_{i=0}^{s-1}{p^{i}}+ (n-s)\cdot p^s)\\
&<1+p^{n-2}+(p-1)(n-(n-2))\cdot p^{n-2}\\
& =
\begin{cases} 1+3\cdot 2^{n-2}< \vert G^\prime \vert & if \qquad p=2; \\
1+5\cdot 3^{n-2}< \vert G^\prime \vert-1 & if\qquad  p=3,
\end{cases}
\end{align*}
which is  contradicts to $t^{L}(KG)=\vert G^{\prime} \vert -p+2$.

Therefore \; $d_{(p^{i}+1)}>0$\; for\; $0\leq i\leq n-2$ and  by
(\ref{e:4}) there exists $\alpha\geq 2$ such that
$d_{(\alpha)}=1$,  and

\begin{align*}
t^{L}(KG)&=
2+(p-1)\sum_{i=0}^{n-2}{p^{i}}+(p-1)(\alpha -1)d_{(\alpha)}\\
&=1+p^{n-1}+(p-1)(\alpha -1)\\
& =
\begin{cases} 2^{n-1}+\alpha & if \qquad p=2; \\
3^{n-1}+2\alpha-1 & if \qquad  p=3.
\end{cases}
\end{align*}
Since $t^{L}(KG)=|G'|-p+2$, it must be $\alpha =p^{n-1}$ and the
proof is done.
\end{proof}

\begin{lemma}\label{L:3}
Let  $K$ be a field with  $\charac(K)=2$ and $G$  a nilpotent
group such that $t^{L}(KG)=\vert G^{\prime}\vert=2^n$. Then one of
the following conditions holds:
\begin{enumerate}
\item[(i)] $cl(G)=2$ and $\gamma_2(G)$  is noncyclic of order $4$;
\item[(ii)] $cl(G)=4$ and $G$ has one of the following properties:
\begin{enumerate}
\item $\gamma_2(G)\cong C_4 \times
C_2$,\quad  $\gamma_3(G)\cong C_2\times C_2$;
\item $\gamma_2(G)\cong  C_2 \times C_2\times
C_2$.
\end{enumerate}
\end{enumerate}
\end{lemma}
\begin{proof} Let $t^{L}(KG)=\vert G^{\prime}\vert=2^n$. Then
either (i) or (ii) of Lemma \ref{L:2} holds. If $n=2$ then by (3)
and (i) of Lemma \ref{L:2} we obtain  that $\frak
D_{(3)}(G)=\gamma_3(G)\cdot \gamma_2(G)^2=\gp{1}$ and the
statement (i) holds.

Assume that $n\geq 3$. By (ii) of Lemma \ref{e:2} we get
$d_{(2^{i}+1)}=1$,\;  $d_{(2^{n-1})}=1$ and  $d_{(j)}=0$, where
\quad $0\leq i\leq n-2$, \quad $j\neq 2^i +1$, \quad $j\neq
2^{n-1}$\quad and $j>1$. The subgroup $H=\frak D_{( 2^{n-1})}(G)$
is central of order $2$ and from (\ref{e:2}) it follows
\begin{eqnarray*}
\frak D_{(m+1)}(G)/H&=&\prod_{(j-1)2^i\geq
m}\gamma_j(G)^{2^i}/H\\
&=&\prod_{(j-1)2^i\geq m}\gamma_j(G/H)^{2^i}=\quad \frak
D_{(m+1)}(G/H).
\end{eqnarray*}
Put $2^{\overline{d}_{(k)}}=[\frak D_{(k)}(G/H):\frak
D_{(k+1)}(G/H)]$ for $k\geq 1$. It is easy to check that
$\overline{d}_{(2^{i}+1)}=1$\quad and \quad
$\overline{d}_{(j)}=0$,  where $0\leq i\leq n-2$, \quad $j\neq 2^i
+1$\quad and $j>1$.

Clearly, $|\gamma_2(G/H)|=2^{n-1}$ and
$t^L(K[G/H])=|\gamma_2(G/H)|+1$. So by Lemma 3 of \cite{BS} and by
Theorem 1 of \cite{BS} the group $\gamma_2(G/H)$ is either a
cyclic $2$-group or $C_2\times C_2$. If $\gamma_2(G/H)$ is  a
cyclic $2$-group, then (a) of Lemma III.7.1 (\cite{H}, p.300)
yields  that   $\gamma_2(G)$ is abelian, so it is isomorphic to
either $C_{2^{n-1}}\times C_2$ or $C_{2^{n}}$. If $\gamma_2(G)$ is
cyclic, then by Theorem 1 of \cite{BS} we get
$t^L(KG)=|G^{\prime}|+1$ and we do not consider this case.

Let $\gamma_2(G/H)=C_2\times C_2$.  It is easy to check that
$|\gamma_2(G)|=8$ and $\gamma_2(G)$ is one of the following
groups: $\mathrm{Q}_{8}$,\; $\mathrm{D}_{8}$,\; $C_4\times C_2$,\;
$C_2\times C_2\times C_2$. It is well known that there do not
exist nilpotent groups $G$ such that  either $\gamma_2(G)\cong
Q_8$ or $\gamma_2(G)\cong D_8$.

Assume that  $\gamma_2(G)=\gp{\;a,b\;\mid\; a^4=b^2=1\;}\cong  C_4
\times C_2$.  Thus by (\ref{e:1})
$$
\frak D_{(3)}(G)=(\frak D_{(2)}(G),G)\cdot \frak
D_{(2)}(G)^2=\gamma_3(G)\cdot \gp{a^2}.
$$
Since   $\vert \frak D_{(2)}(G)/\frak D_{(3)}(G)\vert =2$, only
one of the following cases is  possible:
\[
\begin{aligned}
\gamma_3(G)=\gp{a}&,\qquad
\gamma_3(G)=\gp{ab},\qquad\gamma_3(G)=\gp{a^2,b},\\
&\gamma_3(G)=\gp{a^2b},\qquad \gamma_3(G)=\gp{b}.
\end{aligned}
\]
Now consider separately each  of these:

{\bf Case 1.a}.  Let either $\gamma_3(G)=\gp{a}$ or
$\gamma_3(G)=\gp{ab}$. Since $\gamma_2(G)^2\subset \gamma_3(G)$,
by Theorem III.2.13 (\cite{H}, p.266), we have that $\gamma_k(G)^2
\subseteq \gamma_{k+1}(G)$ for every $k\geq 2$. It follows that
$\gamma_2(G)^2=\gamma_4(G)$. Moreover, $\gamma_{3}(G)^2\subseteq
\gamma_{5}(G)$. Indeed, the elements of the form $(x,y)$, where
$x\in \gamma_2(G)$ and $y\in G$ are generators of $\gamma_3(G)$,
so we have to  prove that $(x,y)^2 \in \gamma_5(G)$. Evidently,
$$
(x^2,y)=(x,y)\cdot (x,y,x)\cdot (x,y)=(x,y)^2\cdot (x,y,x)^{(x,y)}
$$
and  $(x^2,y),\; (x,y,x)^{(x,y)}\in \gamma_5(G)$, so $(x,y)^2\in
\gamma_5(G)$ and $\gamma_3(G)^2\subseteq\gamma_5(G)$. Thus
$\gp{a^2}\subseteq \gp{1}$,  a contradiction.

{\bf Case 1.b}. Let $\gamma_3(G)=\gp{a^2}\times \gp{b}$ and
$cl(G)=3$. Now, let us compute the weak complement of
$\gamma_3(G)$ in $\gamma_2(G)$ (see \cite{BK}, p.34). It is easy
to see that in the notation of \cite{BK} (see p.34)
$$
P=\gamma_2(G),\qquad  H=\gp{a^2,b},\qquad  H\setminus
P^2=\{b,a^2b\}
$$
so $\nu(b)=\nu(a^2b)=2$ and  the weak complement is $A=\gp{a}$.
Since $G$ is of class $3$, by (ii) of Theorem 3.3 (\cite{BK},
p.43) we have
$$
t_L(KG)=t^L(KG)=t(\gamma_2(G))+t(\gamma_2(G)/\gp{a})=7\not=
|G^\prime|.
$$

{\bf Case 1.c}. Let either $\gamma_3(G)=\gp{b}$ or
$\gamma_3(G)=\gp{a^2b}$. Clearly $cl(G)=3$.   According to the
notation of \cite{BK} (see p.34) we have that  $P=\gamma_2(G)$ and
either $H=\gp{b}$ and $H\setminus P^2=\{b\}$ or $H=\gp{a^2b}$ and
$H\setminus P^2=\{a^2b\}$. It follows that $\nu(b)=\nu(a^2b)=2$
in  both cases and the weak complement is $A=\gp{a}$. As  in the
case 1.b we have $t^L(KG)=7\not= |G^\prime|$, a contradiction.

Now,  let $\gamma_2(G)\cong C_2\times C_2\times C_2$. If
$\gamma_3(G)\cong C_2$ then by (\ref{e:1})
$$
\frak D_{(2)}(G)=\gamma_2(G),\qquad \frak
D_{(3)}(G)=\gamma_3(G),\qquad \frak D_{(4)}(G)=\gp{1}
$$
and $d_{(2)}=2$,\quad  $d_{(3)}=1$, which  contradicts (ii) of
Lemma \ref{L:2}.

If $\gamma_3(G)\cong C_2\times C_2$ and $cl(G)=3$, then by
(\ref{e:1}) we have
$$
\frak D_{(2)}(G)=\gamma_2(G),\qquad \frak
D_{(3)}(G)=\gamma_3(G),\qquad \frak D_{(4)}(G)=\gp{1},
$$
also  a contradiction, because $d_{(2)}=1$\quad and \quad
$d_{(3)}=2$.

Finally, let $\gamma_2(G)=\gp{a,b\mid a^{2^{n-1}}=b^2=1}\cong
C_{2^{n-1}}\times C_2$ with $n\geq 4$. According to (1), $\frak
D_{(2)}(G)=\gamma_2(G)$\; and \; $\frak
D_{(3)}(G)=\gamma_3(G)\cdot\gp{a^2}$ hold. Since   $\vert \frak
D_{(2)}(G)/\frak D_{(3)}(G)\vert =2$,  we obtain one of the
following cases:
\[
\begin{aligned}
\gamma_3(G)=\gp{a}&,\qquad\qquad \gamma_3(G)=\gp{ab},\qquad\qquad\gamma_3(G)=\gp{b},\\
& \gamma_3(G)=\gp{a^{2^j}b},\qquad \gamma_3(G)=\gp{a^{2^j},
b},\\
\end{aligned}
\]
where $1\leq j\leq n-2$.

We consider each of these:

{\bf Case 2.a}. Let either $\gamma_3(G)=\gp{a}$ or
$\gamma_3(G)=\gp{ab}$ or $\gamma_3(G)=\gp{b}$.    Using the
arguments of the cases 1.a and 1.c above, it is easy to verify
that we obtain a contradiction.

{\bf Case 2.b}.  Let $\gamma_3(G)=\gp{a^{2^{j}}b}$ with  $1\leq
j\leq n-2$.  Then by (\ref{e:1}) and by (ii) of Lemma \ref{L:2} we
get
\[
\begin{aligned}
\frak D_{(2)}(G)&=\gamma_{2}(G);\qquad \qquad \frak D_{(3)}(G)=\gp{a^{2}}\times \gp{b};\\
\frak D_{(4)}(G)&=(\gp{a^{2}}\times \gp{b},G)\cdot\gp{a^{4}}=
\begin{cases}  \gp{a^{2}};\\ \gp{a^{2}b};\\ \gp{a^{4}}\times
\gp{b}.
\end{cases}
\end{aligned}
\]
Suppose that $\frak D_{(4)}(G)=\gp{a^{2}}$. By (ii) of Lemma
\ref{e:3}
$$
\frak D_{(5)}(G)=(\gp{a^{2}},G)\cdot \frak
D_{(3)}(G)^2=(\gp{a^{2}},G)\cdot \gp{a^{4}}=\gp{a^{2}}.
$$
The last equality forces $(\gp{a^{2}},G)=\gp{a^{2}}$ and so $\frak
D_{(k)}(G)=\gp{a^{2}}$ for each $k\geq 5$, which is impossible.

Now let  $\frak D_{(4)}(G)=\gp{a^{2}b}$. As above,
$$
\frak D_{(4)}(G)=\frak D_{(5)}(G)=(\gp{a^{2}b},G)\cdot
\gp{a^{4}}=\gp{a^{2}b},
$$
and  we get $(\gp{a^{2}b},G)=\gp{a^{2}b}$, which is not possible
either.

Finally, let  $\frak D_{(4)}(G)=\gp{a^{4}}\times \gp{b}$. Suppose
that there exists $k \leq 2^{n-2}+1$, such that $\frak D_{(k)}(G)$
is cyclic. Using the same arguments as above,  we obtain that
$\frak D_{(m)}(G)\not=\gp{1}$ for each $m$, which is impossible.

So $\frak D_{(2^{n-2}+1)}(G)=\gp{a^{2^{n-2}}}\times \gp{b}$ and by
(\ref{e:1}) and by (ii) of Lemma \ref{L:2}
\[
\begin{aligned}
\frak D_{(2^{n-2}+2)}(G)&=(\frak D_{(2^{n-2}+1)}(G), G )\cdot
\frak D_{(2^{n-3}+2)}(G)^2\\
&=(\frak D_{(2^{n-2}+1)}(G), G )=\gp{\omega\mid \omega^2=1};\\
\frak D_{(2^{n-2}+3)}(G)&=(\frak
D_{(2^{n-2}+2)}(G),G)=(\gp{\omega},G)=\gp{\omega},
\end{aligned}
\]
which is not possible either.

{\bf Case 2.c}. Let $\gamma_3(G)=\gp{a^{2^{j}}}\times \gp{b}$ with
$1\leq j\leq n-2$. It is easy to check that this case is similar
to the last subcase of the case 2.b.

So, the proof is complete.\end{proof}

\begin{lemma}\label{L:4}
Let  $K$ be a field with  $\charac(K)=3$ and $G$  a nilpotent
group  such that $\vert G^{\prime}\vert=3^n$  and $t^{L}(KG)=\vert
G^{\prime}\vert-1$. Then
 $cl(G)=3$, $\gamma_2(G)\cong C_3\times C_3$ and $\gamma_3(G)\cong C_3$.
\end{lemma}
\begin{proof} Let $t^{L}(KG)=3^n-1$. Then
either (iii) or (iv) of Lemma \ref{L:2} holds.  By (iv) of Lemma
\ref{L:2}  it yields
$$
d_{(3^{i}+1)}=1,\qquad  d_{(3^{n-1})}=1, \qquad d_{(j)}=0,
$$
where \quad $0\leq i\leq n-2$, \quad $j\neq 3^i +1$, \quad $j\neq
3^{n-1}$\quad and $j>1$.

The subgroup $H=\frak D_{( 3^{n-1})}(G)$ is  central of order $3$
and from (\ref{e:2}), as we already proved at the beginning of the
proof of Lemma \ref{L:3}, we have
\begin{eqnarray*}
\frak D_{(m+1)}(G)/H=\frak D_{(m+1)}(G/H).
\end{eqnarray*}
It follows that\quad  $\overline{d}_{(3^{i}+1)}=1$\quad for $0\leq
i\leq n-2$, \quad $\overline{d}_{(j)}=0$\quad for $j\neq 3^i
+1$\quad and $j>1$,  where \quad $3^{\overline{d}_{(k)}}=[\frak
D_{(k)}(G/H):\frak D_{(k+1)}(G/H)]$.

Clearly, $|\gamma_2(G/H)|=3^{n-1}$ and
$t^L(K[G/H])=|\gamma_2(G/H)|+1$. By Lemma 3 of  \cite{BS} and by
Theorem 1 of \cite{BS} we have that $\gamma_2(G/H)$ is  a cyclic
$3$-group. By (a) of Lemma III.7.1 (\cite{H}, p.300) the group
$\gamma_2(G)$ is abelian, so it is isomorphic to either
$C_{3^{n-1}}\times C_3$ or $C_{3^{n}}$. By Theorem 1 of \cite{BS},
in the last case $KG$ has upper and lower maximal Lie nilpotency
index.

Therefore we can assume that $\gamma_2(G)\cong C_{3^{n-1}}\times
C_3$ and $n\geq 3$. Since $exp(\gamma_2(G))=3^{n-1}$\; and
$3^{n-2}|\; 2\cdot 3^{n-2}$, \; by (iii) of Lemma \ref{L:2} and by
(ii) of Lemma \ref{sha} we obtain that $\frak
D_{(3^{n-2}+1)}(G)=\gp{1}$ and so $\frak D_{(3^{n-1})}(G)=\gp{1}$
in contradiction with (iii) of Lemma \ref{L:2}. It follows that
$\gamma_2(G)\cong C_{3}\times C_3$.

Suppose  that $\gamma_2(G)\subseteq \zeta(G)$. By (\ref{e:1}) we
get $\frak D_{(3)}(G)=\gp{1}$ and so $d_{(2)}=2$, which is in
contradiction with (iii) of Lemma \ref{L:2}.

Finally, the case (iii) of Lemma \ref{L:2} belongs  to the
previous case. The proof is done.
\end{proof}
 Proof of the Theorem  follows from   Lemma 3 and Lemma 4. The
converse   is trivial.

\end{document}